\documentclass[11pt, a4paper, oneside, reqno]{amsart} 
\usepackage{amssymb, amsmath}% makeidx}

\pdfoutput=1
\usepackage{graphicx}
\usepackage{amssymb, amsmath, amsthm}
 \usepackage[bookmarks=false]{hyperref}

\newcommand*{\fplus}{\genfrac{}{}{0pt}{}{}{+}}
\newcommand*{\fdots}{\genfrac{}{}{0pt}{}{}{\cdots}}

\makeatletter
\renewcommand{\@biblabel}[1]{#1.} 
\makeatother

\newtheorem{Theorem}[equation]{Theorem}

\newtheorem*{rem}{Remark} %\numberwithin{equation}{section}

\newtheorem*{rems}{Remarks} %\numberwithin{equation}{section}

\newtheorem*{sol}{Solution}

\newtheorem*{nt}{Notes}

\makeatletter
\renewcommand{\@biblabel}[1]{#1.}
\makeatother

\newcommand{\qrfac}[2]{{\left({#1}; q\right)_{#2}}}

\numberwithin{equation}{section}

\begin{document} 
\title[Convergents of  a continued fraction of Ramanujan]{A formula for the convergents of  a continued fraction of Ramanujan}

\author[Gaurav Bhatnagar]{Gaurav Bhatnagar$^*$
}
\address{Indian Statistical Institute, Delhi Center, 7 S.~J.~S.~Sansanwal Marg, Delhi. }
\curraddr{Fakult\"at f\"ur Mathematik,  Universit\"at Wien \\
Oskar-Morgenstern-Platz 1, 1090 Wien, Austria.}
\email{bhatnagarg@gmail.com}
\thanks{$^*$Research is supported in part by the Austrian Science Fund (FWF) grant F 5008-N15. 
}

\author{Michael D. Hirschhorn}
\address{School of Mathematics and Statistics, UNSW, Australia 2052}
\email{m.hirschhorn@unsw.edu.au}

\date{\today}

%\keywords{Rogers-Ramanujan Continued Fraction,  Ramanujan, the Lost Notebook}
%\subjclass[2010]{Primary 33D15;  Secondary 11A55,  30B70}

\begin{abstract}
In Entry 16, Chapter 16 of his notebooks, Ramanujan himself gave a formula for the convergents of the famous Rogers-Ramanujan continued fraction. We provide a similar formula for the convergents of a more general continued fraction, namely Entry 15 of Chapter 16. \\
{\bf Keywords}: Rogers-Ramanujan Continued Fraction,  Ramanujan, the Lost Notebook.
\end{abstract}

\maketitle

\section{The Theorem}
The $q$-rising factorial $(a;q)_k$ is defined as $(a;q)_0:=1$; and when $k>0$, as the product of $k$ terms:
$$(a;q)_k=(1-a)(1-aq)\cdots (1-aq^{k-1}).$$
Ramanujan wrote a ratio of two series as a continued fraction:
\begin{equation}\label{cor-entry15-a}
\frac{\displaystyle\sum_{k=0}^{\infty} \frac{q^{k^2}}{\qrfac{q}{k}}\lambda^k}
{\displaystyle\sum_{k=0}^{\infty} \frac{q^{k^2+k}}{\qrfac{q}{k}}\lambda^k}
=
1+\frac{\lambda q}{1}\fplus\frac{\lambda q^2}{1}\fplus\frac{\lambda q^3}{1}\fplus\fdots.
\end{equation}
This is known as the Rogers-Ramanujan continued fraction (the two series on the left are the sum sides of the famous Rogers-Ramanujan identities!) and appears as a corollary to Entry 15 in Chapter 16 of Ramanujan's Notebooks 
(see Berndt \cite{berndt-notebooks-3}).   Next, in Entry 16, Ramanujan provides a formula for the $n$th convergents of this continued fraction:
For each positive integer $n$, let
\begin{equation*}\label{mu-entry16}
\mu =\mu_n(\lambda, q) =
\sum_{k=0}^{[(n+1)/2]} \frac{q^{k^2}\lambda^k}{\qrfac{q}{k}}
\frac{\qrfac{q}{n-k+1}}{\qrfac{q}{n-2k+1}}
\end{equation*}
and
\begin{equation*}\label{nu-entry16}
\nu =\nu_n(\lambda, q) =
\sum_{k=0}^{[n/2]} \frac{q^{k^2+k}\lambda^k}{\qrfac{q}{k}}
\frac{\qrfac{q}{n-k}}{\qrfac{q}{n-2k}}.
\end{equation*}
Then,
\begin{equation}
\frac{\mu}{\nu} =
%=1+ \frac{\lambda q}{(1+bq)}\frac{\mu_n(2)}{\mu_n(1)}=
1+ \frac{\lambda q}{1}\fplus
\frac{\lambda q^2}{1}\fplus\frac{\lambda q^3}{1}
\fplus\fdots \fplus \frac{\lambda q^n}{1}
 .\label{entry16}
\end{equation}

Ramanujan's generalization of the Rogers-Ramanujan continued fraction is given by Entry 15 of Chapter 16 of 
\cite{berndt-notebooks-3}: For $|q|<1$, 
%\begin{align}
 \begin{equation}
 \frac{\displaystyle
\sum_{k=0}^{\infty}\frac{q^{k^2}}{\qrfac{q}{k}\qrfac{-bq}{k}}\lambda^k}
{\displaystyle\sum_{k=0}^{\infty}\frac{q^{k^2+k}}{\qrfac{q}{k}\qrfac{-bq}{k}} \lambda^k}
%\label{g-frac-sums-o}
=1+\frac{\lambda q}{1+bq}\fplus\frac{\lambda q^2}{1+bq^2}\fplus\frac{\lambda q^3}{1+bq^3}\fplus\fdots 
\label{g-cfrac2-o}
\end{equation}
%\end{align}
where we have renamed some symbols in order to  fit the notation used by
Andrews and Berndt~\cite[Entry 6.3.1 (ii)]{LN1}. 

However, Ramanujan did not provide a formula for the convergents of this continued fraction. Our objective here is to provide just such a formula, which is a natural extension of Ramanujan's own formula \eqref{entry16}.  
\begin{Theorem}\label{mu1-finite-thm}
Let
\begin{equation}\label{mu1-finite}
g_n(s):=
\sum_{k=0}^{\infty} \frac{q^{k^2+sk}\lambda^k}{\qrfac{q}{k}\qrfac{-bq^s}{k}}
\frac{\qrfac{q}{n-k-s+1}}{\qrfac{q}{n-2k-s+1}}
\frac{\qrfac{-bq}{n-k}}{\qrfac{-bq}{n}}.
\end{equation}
Then, for $n=1, 2, 3, \dots$, we have
\begin{equation}
(1+b)\frac{g_n(0)}{g_n(1)} =
%=1+ \frac{\lambda q}{(1+bq)}\frac{\mu_n(2)}{\mu_n(1)}=
1+b + \frac{\lambda q}{1+bq}\fplus
\frac{\lambda q^2}{1+bq^2}\fplus\frac{\lambda q^3}{1+bq^3}
\fplus\fdots \fplus \frac{\lambda q^n}{1+bq^n}
 .\label{entry16-gen1-a}
\end{equation}
\end{Theorem}

Observe that when $b=0$, then \eqref{entry16-gen1-a} reduces to \eqref{entry16}. In this case, $g_n(0)$ and $g_n(1)$ reduce to $\mu$ and $\nu$ respectively. 

To obtain a formula for the convergents of \eqref{g-cfrac2-o}, consider:
\begin{equation*}
(1+b)\frac{g_n(0)}{g_n(1)} -b
 = 1+  \frac{\lambda q}{1+bq}\fplus
\frac{\lambda q^2}{1+bq^2}\fplus\frac{\lambda q^3}{1+bq^3}
\fplus\fdots \fplus \frac{\lambda q^n}{1+bq^n}
 .%\label{g-cfrac2-finite-a}
\end{equation*}

\section{A Proof}
Before heading into the proof of \eqref{entry16-gen1-a}, we need the definition of $q$-rising factorials, when $k$ is not a non-negative integer.  For that we need the infinite $q$-rising factorial
$$\qrfac{a}{\infty}:=\prod_{j=0}^{\infty} (1-aq^j), {\text{ for $|q|<1$}}.$$
When $k$ is not a positive integer, one can define
$$(a;q)_k=\frac{\qrfac{a}{\infty}}{\qrfac{aq^k}{\infty}}.$$
Observe that this definition implies
$$\frac{1}{\qrfac{q}{m}} = 0 \text{ when } m=-1, -2, -3, \dots .$$ 

To prove \eqref{entry16-gen1-a}, we use the approach used by Euler~\cite{eulerE616} (as explained by Bhatnagar \cite{GB2014}). We use the elementary identity:
\begin{equation}\label{div-1step}
\frac{N}{D}=1+\frac{N-D}{D}
\end{equation}
to \lq divide' two series. 

\begin{proof}
Observe first that in the sum $g_n(s)$ in \eqref{mu1-finite},  the index $k$ goes from $0$ to $ \lfloor \frac{n-s+1}{2}\rfloor$.
Further, observe that
\begin{equation}\label{endpoints}
g_n(n)=1=g_n(n+1),
\end{equation} 
since only the terms corresponding to the index $k=0$ survive.  

We will show
\begin{equation}\label{gn-recursion}
(1+bq^s)\frac{g_n(s)}{g_n(s+1)}=
1+bq^s+\frac{\lambda q^{s+1}}{(1+bq^{s+1}){\cfrac{g_n(s+1)}{g_n(s+2)}}},
\end{equation}
for $s=0,1, 2, 3, \dots, n-1$.
The formula \eqref{entry16-gen1-a} follows immediately by iterating \eqref{gn-recursion} $n$ times.

To prove \eqref{gn-recursion}, we 
use \eqref{div-1step} to find that
\begin{equation}\label{proof-step1}
(1+bq^s) \frac{g_n(s)}{g_n(s+1)}=(1+bq^s)\left(1+\frac{g_n(s)-g_n(s+1)}{g_n(s+1)}\right).
\end{equation}
Consider the difference of sums $g_n(s)-g_n(s+1)$.
\begin{align*}
g_n(s) -  & g_n(s+1) = 
\displaystyle\sum_{k=0}^{\infty}  \Bigg(
\frac{q^{k^2+sk}\lambda^k}{\qrfac{q}{k}\qrfac{-bq^s}{k+1}}
\frac{\qrfac{q}{n-k-s}}{\qrfac{q}{n-2k-s+1}}
\frac{\qrfac{-bq}{n-k}}{\qrfac{-bq}{n}} \notag
\\
&\times  \left[(1+bq^{s+k}) (1-q^{n-k-s+1})-  (1+bq^s)(1-q^{n-2k-s+1})q^k\right]\Bigg) \notag\\
&= \displaystyle\sum_{k=0}^{\infty}  \Bigg(
\frac{q^{k^2+sk}\lambda^k}{\qrfac{q}{k}\qrfac{-bq^s}{k+1}}
\frac{\qrfac{q}{n-k-s}}{\qrfac{q}{n-2k-s+1}}
\frac{\qrfac{-bq}{n-k}}{\qrfac{-bq}{n}} \notag
\\
&\times  \left[(1+bq^{n-k+1}) (1-q^{k})\right]\Bigg) \notag \\
&= \displaystyle\sum_{k=1}^{\infty}  
\frac{q^{k^2+sk}\lambda^k}{\qrfac{q}{k-1}\qrfac{-bq^s}{k+1}}
\frac{\qrfac{q}{n-k-s}}{\qrfac{q}{n-2k-s+1}}
\frac{\qrfac{-bq}{n-k+1}}{\qrfac{-bq}{n}} \notag \\
&= \displaystyle\sum_{k=0}^{\infty}  
\frac{q^{(k+1)^2+s(k+1)}\lambda^{k+1}}{\qrfac{q}{k}\qrfac{-bq^s}{k+2}}
\frac{\qrfac{q}{n-k-s-1}}{\qrfac{q}{n-2k-s-1}}
\frac{\qrfac{-bq}{n-k}}{\qrfac{-bq}{n}} \notag \\
&=\frac{\lambda q^{s+1}}{(1+bq^s)(1+bq^{s+1})}g_n(s+2).% \label{proof-step2}
\end{align*}
Now by substituting  in \eqref{proof-step1}, we immediately obtain \eqref{gn-recursion}, and our proof is complete.
\end{proof}

Notice that Ramanujan's continued fraction \eqref{g-cfrac2-o} is an immediate corollary of our formula. Take the limits $n \to \infty$ in \eqref{entry16-gen1-a}. The continued fraction in \eqref{g-cfrac2-o} converges if its convergents converge, and we have  found the convergents to be $(1+b)g_n(0)/g_n(1)-1$ which converges when $|q|<1$. (Getting the left hand side of \eqref{g-cfrac2-o} from this is a pleasant exercise. Try it!)

\section{Some connections}
Ramanujan's Entry 16 (equation \eqref{entry16}) was rediscovered by P. Kesava Menon \cite{pkmenon65}, and Hirschhorn \cite{mdh72}. Another formula for the convergents of an even more general continued fraction of Ramanujan has been given by Hirschhorn \cite[eq. (1)]{mdh74}. On taking $a=0$ in Hirschhorn's formula, we obtain a formula (different from \eqref{entry16-gen1-a}) for the convergents of \eqref{g-cfrac2-o}, where the numerator and denominator are double sums. Similar results appear in Hirschhorn \cite{mdh80} and \cite{mdh92}.

There is also a sequence of orthogonal polynomials due to Al-Salam and Ismail \cite{al-salam-ismail1983}, namely
\begin{equation*}
U_n(x;a,b)= \sum_{k\geq 0} \frac{\qrfac{-a}{n-k}\qrfac{q}{n-k}}{\qrfac{-a}{k}\qrfac{q}{k}\qrfac{q}{n-2k}}
x^{n-2k} (-b)^k q^{k(k-1)}.
\end{equation*}
This is similar to our $g_n(s)$. Indeed, we can see  that
\begin{equation*}
g_n(1)=\frac{1}{\qrfac{-bq}n}U_n(1;bq,-\lambda q^2).
\end{equation*}
Al-Salam and Ismail~\cite{al-salam-ismail1983} have not explicitly stated \eqref{entry16-gen1-a}, but perhaps our formula can be extracted from their work. 

%There are many more $q$-continued fractions of Ramanujan that can be proved by using Euler's approach, and %have recurrence relations such \eqref{gn-recursion}. Perhaps the connection with orthogonal polynomials can %help in finding formulas for their convergents. 

\bibliographystyle{amsplain}

\end{document}